\documentclass[10pt]{amsart}
\usepackage{latexsym, mathrsfs, color, multirow,bbm,mathtools,amsmath,amssymb,tikz-cd}
\usepackage{amsmath}
\usepackage{graphics}
\usepackage{subcaption}
\usepackage{xparse}
\usepackage{faktor}
\usepackage{soul}
\usepackage{multirow}
\input{xy}
\xyoption{all}

\usepackage[latin1]{inputenc}
\usepackage{tikz}

\usepackage[colorlinks=true, pdfstartview=FitV,
 linkcolor=blue,citecolor=black,urlcolor=black]{hyperref}
 
\usepackage{tikz-cd}
\usepackage{tkz-euclide}

\numberwithin{equation}{section}
\theoremstyle{plain}

\newcounter{problems}
\theoremstyle{plain}
\newtheorem{theorem}{Theorem}[section]

\newtheorem{corollary}[theorem]{Corollary}

\theoremstyle{definition}
\newtheorem{example}[theorem]{Example}

\newtheorem*{remark}{Remark}
\theoremstyle{definition}
\newtheorem{Def}[equation]{Definition}

\newenvironment{red}{\relax\color{red}}{\relax}
\newenvironment{blue}{\relax\color{blue}}{\hspace*{.5ex}\relax}

\newcommand{\ber}{\begin{red}}
\newcommand{\er}{\end{red}}
\newcommand{\beb}{\begin{blue}}
\newcommand{\eb}{\end{blue}}

\newcommand{\bw}{{\boldsymbol{w}}}

\usepackage[nodisplayskipstretch]{setspace}  

\begin{document}

\title{An unexpected property of $\mathbf{g}$-vectors for rank 3 mutation-cyclic quivers}

\author{Jihyun Lee}
\address{Department of Mathematics, University of Alabama,
	Tuscaloosa, AL 35487, U.S.A. 
	}
\email{jlee291@crimson.ua.edu}

\author{Kyungyong Lee}
\address{Department of Mathematics, University of Alabama,
	Tuscaloosa, AL 35487, U.S.A. 
	and Korea Institute for Advanced Study, Seoul 02455, Republic of Korea}
\email{kyungyong.lee@ua.edu; klee1@kias.re.kr}

\thanks{The authors were supported by the University of Alabama, Korea Institute for Advanced Study, and the NSF grant DMS 2042786 and DMS 2302620.}

\date{\today}


\begin{abstract}
Let $Q$ be a rank 3 mutation-cyclic quiver. 
It is known that every $\mathbf{c}$-vector of $Q$ is a solution to a quadratic equation of the form 
$$\sum_{i=1}^3 x_i^2 + \sum_{1\leq i<j\leq 3} \pm q_{ij} x_i x_j =1,$$
where $q_{ij}$ is the number of arrows between the vertices $i$ and $j$ in $Q$. A similar property holds for $\mathbf{c}$-vectors of any acyclic quiver. 
In this paper, we show that $\mathbf{g}$-vectors of $Q$ enjoy an unexpected property. More precisely, every $\mathbf{g}$-vector of $Q$ is a solution to a quadratic equation of the form 
$$\sum_{i=1}^3 x_i^2 + \sum_{1\leq i<j\leq 3} p_{ij} x_i x_j =1,$$
where $p_{ij}$ is the number of arrows between the vertices $i$ and $j$ in another quiver $P$  obtained by mutating $Q$.
\end{abstract}
\maketitle

\section{Introduction}

In the context of cluster algebras, $\mathbf{g}$-vectors, along with $\mathbf{c}$-vectors, serve as fundamental objects. The study of $\mathbf{g}$-vectors is crucial for understanding the properties of cluster variables.

When a quiver $Q$ is acyclic, it is known that the positive $\mathbf{c}$-vectors of $Q$ are real Schur roots, that is, the dimension vectors of indecomposable rigid modules over $Q$ \cite{chavez2015c, hubery2016categorification, igusa2010exceptional, seven2015cluster, speyer2013acyclic}. In particular, every $\mathbf{c}$-vector of $Q$ is a real root, hence  is a solution to a quadratic equation of the form 
$$\sum_{i=1}^n x_i^2 + \sum_{1\leq i<j\leq n} \pm q_{ij} x_i x_j =1,$$
where $n$ is the number of vertices of $Q$ and $q_{ij}$ is the number of arrows between the vertices $i$ and $j$ in $Q$. It is also shown \cite{EJLN2024, seven2024cluster} that the same property holds for $\mathbf{c}$-vectors of rank 3 mutation-cyclic quivers. 

 For quivers of rank at most 2, explicit formulas for $\mathbf{g}$-vectors are well-understood. This paper focuses on the $\mathbf{g}$-vectors of rank 3 mutation-cyclic quivers.

Let $Q$ be a rank 3 mutation-cyclic quiver. 
Our main theorem states that for any reduced mutation sequence $\bw$ from $Q$ whose first mutation occurs at vertex $k \in \{1,2,3\}$, each column vector $\mathbf{g}_i^{\bw}$ in the $G^\bw$-matrix, referred to as a $\mathbf{g}$-vector, satisfies the quadratic equation $$q(\mathbf{g}_{i}^{\bw}) = x_1^2 + x_2^2 + x_3^2 + |b_{12}^{[k]}|x_1x_2 + |b_{23}^{[k]}|x_2x_3 + |b_{31}^{[k]}|x_3x_1 = 1,$$ where the coefficients are determined by the number of arrows between vertices $i$ and $j$ in the quiver $Q^{[k]}$ obtained by mutating $Q$ at the vertex $k$.

\noindent\emph{Acknowledgements.} We would like to thank Tucker Ervin, Blake Jackson, Jae-Hoon Kwon, Kyu-Hwan Lee, and Ahmet Seven for the correspondences.

\section{Preliminaries} \label{sec-preliminaries}

Let $n$ be a positive integer. To a quiver $Q$ with vertices labeled $1,2,\cdots,n,$ we can associate an $n\times n$ skew-symmetric matrix, i.e. the matrix is in correspondence with a quiver $Q.$ 

\begin{Def}\label{Bmatrix}
Let $B = [b_{ij}]$ be an $n\times n$ skew-symmetric matrix. This $B$ will be called an \emph{exchange matrix}. This matrix determines a quiver $Q$ as follows: $b_{ij} > 0$ if and only if $Q$ has $b_{ij}$ arrows from vertex $i$ to vertex $j$. Sometimes we write $B=B(Q)$.
\end{Def}

Then, we extend the notion of mutation from quivers to (exchange) matrices.

\begin{Def}\emph{(Matrix mutation)} \label{Matrixmutation}\\
We review the definition of matrix mutations, following \cite{fomin2002cluster}. 
Now assume that $M=[m_{ij}]$ is an $2n \times n$ matrix of integers. Let $\mathcal I:= \{ 1, 2, \dots, n \}$ be the set of indices. For $\bw=[i_1, i_2, \dots , i_t]$, $i_j \in \mathcal I$, we define the matrix $M^\bw=[m_{ij}^\bw]$ inductively: the initial matrix is $M$ for $\bw=[\,]$, and assuming we have $M^\bw$, define the matrix $M^{\bw[k]}=[m_{ij}^{\bw[k]}]$ for $k \in \mathcal I$ with $\bw[k]:=[i_i, i_2, \dots , i_t, k]$
by  
\begin{equation} \label{eqn-mmuu} m_{ij}^{\bw[k]} = \begin{cases} -m_{ij}^\bw & \text{if  $i=k$ or $j=k$}, \\ m_{ij}^\bw + \mathrm{sgn}(m_{ik}^\bw) \, \max(m_{ik}^\bw m_{kj}^\bw,0) & \text{otherwise}, \end{cases}
\end{equation}
where  $\mathrm{sgn}(a) \in \{1, 0, -1\}$ is the signature of $a$. The matrix $M^{\bw[k]}$ is called the {\em mutation of $M^\bw$} at index (or label) $k$, $\bw$ and $\bw[k]$ are called {\em mutation sequences}, and $n$ is the {\em rank}.
If the mutation sequence $[i_1, i_2, \dots , i_t]$ satisfies $i_j\neq i_{j+1}$ for all $j\in\{1,...,n-1\}$, then it is said to be \emph{reduced}.
\end{Def}

\begin{Def} \emph{($C$-matrix, $\mathbf{c}$-vector and $G$-matrix, $\mathbf{g}$-vector)}\\
We also recall the definitions of the $\mathbf{c}$-vector and $\mathbf{g}$-vector as introduced in \cite{fomin2002cluster, fomin2007cluster}. Let $B$ be an $n\times n$ initial exchange matrix and $C$ be an initial $C$-matrix. Consider the $2n\times n$ matrix $M=\begin{pmatrix}B\\
C\end{pmatrix}$ and
 a mutation sequence $\bw=[i_1, \dots ,i_t]$. After the mutation sequence $\bw$, i.e., at the indices $i_1, \dots , i_t$ consecutively, we obtain $\begin{pmatrix}B^\bw \\
 C^\bw\end{pmatrix}$. Write their entries as $B^\bw= \begin{bmatrix} b_{ij}^\bw \end{bmatrix}, C^\bw= \begin{bmatrix} c_{ij}^\bw \end{bmatrix} = \begin{bmatrix}\mathbf{c}_1^\bw & \cdots & \mathbf{c}_n^\bw \end{bmatrix},$ where $\mathbf{c}_i^\bw$ are the column vectors. The resulting column vectors $\mathbf{c}_i^\bw$ are called \emph{$\mathbf{c}$-vectors} of $B$ for $\bw,$ and the matrix $C^\bw$ is called a \emph{$C$-matrix} of $B$ for $\bw.$ And $G^\bw = ((C^\bw)^{-1})^T$ matrix is called a \emph{$G$-matrix} of $B$ for $\bw.$ From the matrix $G^\bw= \begin{bmatrix} g_{ij}^\bw \end{bmatrix} = \begin{bmatrix}\mathbf{g}_1^\bw & \cdots & \mathbf{g}_n^\bw \end{bmatrix},$ each column vector $\mathbf{g}_i^\bw$ is called a \emph{$\mathbf{g}$-vector} of $B$ for any $\bw.$
\end{Def}
Although this definition of the $G$-matrix is clearly established, it requires computing the inverse of a matrix, which can be cumbersome. Additionally, for proving our main theorem, this approach is less convenient, as it complicates the recursive relationships between $G$-matrices during mutations. Therefore, following the notation developed in \cite{nakanishi2012tropical}, we adopt a modified approach to defining $G$-matrix mutation that better suits our needs.

Following \cite{nakanishi2012tropical}, we first extend the notation $[b]_+=\text{max} (b,0)$ to matrices, writing $[B]_+$ for the matrix obtained from $B$ by applying the operation $b\mapsto[b]_+$ to all entries of $B.$ For matrix index $l,$ we denote by $B^{\bullet l}$ the matrix obtained from $B$ by replacing all entries outside of the $l$-th column with zeros; the matrix $B^{l \bullet}$ is defined similarly using the $l$-th row instead of the column.

Note that the operations $B\mapsto[B]_+$ and $B\mapsto B^{\bullet l}$ commute with each other, making the notation $[B]_+^{\bullet l}$(and $[B]_+^{l \bullet}$) straightforward.
Using this formalism, we can rephrase the sign-coherence condition \cite{GHKK} for a $C$-matrix as follows :
\begin{equation}\label{Csigncoherence}
\text{For every $j,$ there exists the sign}\quad \varepsilon_j(C) = \pm 1 \quad \text{such that} \quad [-\varepsilon_j(C)C]_+^{\bullet j} = 0.  
\end{equation}

For each mutation sequence $\bw=[i_1, i_2, \dots , i_t]$, we then inductively define the $C,G$-matrix of $B$ as follows. The following definition provides a recursive method to find the $C$ and $G$-matrices by using the previous $C$ and $G$-matrices.
\begin{Def} \emph{(C,G-matrix mutation)} \label{CGmatrix}\\
Let $C$ and $G$ be the initial $C$ and $G$-matrices, both of which are $n \times n$ identity matrices. Given a mutation sequence $\bw=[i_1, i_2, \dots , i_t]$, suppose that $C^\bw$and $G^\bw$ are defined. We aim to define $C^{\bw[k]}$ and $G^{\bw[k]}$ such that $\mu_k(B^\bw, C^\bw, G^\bw) = (B^{\bw[k]}, C^{\bw[k]} , G^{\bw[k]}).$

Then, under the assumption that $C$ satisfies (\ref{Csigncoherence}), we have
\begin{equation} \label{CG-mutation}
C^{\bw[k]}=C^{\bw} (J_k + [\varepsilon_k(C^{\bw})B^{\bw}]_+^{k \bullet}), \quad
G^{\bw[k]}=G^{\bw} (J_k + [-\varepsilon_k(C^{\bw})B^{\bw}]_+^{\bullet k}),
\end{equation}
where we use the notation $J_k$ to denote the diagonal matrix obtained from the identity matrix by replacing the $(k,k)$-entry with $-1.$
\end{Def}

\subsection{Examples of Definition \ref{CGmatrix}}\label{example2.1}

As examples of how to calculate $C$ and $G$-matrices in Definition \ref{CGmatrix}, consider the following quiver $Q$ as an initial quiver:
\begin{equation}\label{rank3initialquiver}
Q=\vcenter{\hbox{
\begin{tikzpicture}[->,>=stealth',shorten >=1pt,auto,node distance=1.5cm,thick,main node/.style={circle,draw}]
  \node[main node] (2) {2};
  \node[main node] (1) [below left of=2] {1};
  \node[main node] (3) [below right of=2] {3};

  \path
    (2) edge [left] node[auto, midway, above] {3} (1)  
    (3) edge [left] node[auto, midway, above] {3} (2)
    (1) edge [right] node[auto, midway, above] {3} (3);
\end{tikzpicture}}},
\hspace{0.5cm}
B(Q) = \left[\begin{smallmatrix}
0 & -3 & 3 \\ 3 & 0 & -3 \\ -3 & 3 & 0
\end{smallmatrix}\right]
\end{equation}

Then the initial framed quiver $F$ corresponding to the given quiver $Q$ with an initial $C$-matrix, which is the identity matrix, is represented by the matrix $M:$
$$M={\scriptsize \left [ \begin{array}{rrr} 0&-3&3\\3&0&-3 \\-3&3&0 \\ \hline 1&0&0\\0&1&0 \\0&0&1\end{array} \right],}\hspace{0.2cm} C = \left[\begin{smallmatrix}
1 & 0 & 0 \\ 0 & 1 & 0 \\ 0 & 0 & 1
\end{smallmatrix}\right].$$
Obviously, the initial $G$-matrix is also the identity matrix.

We will consider the following mutation sequences : $\bw=[2]$, $[2,3]$, or $[2,3,1]$. 

\noindent\textbf{i)} $\mu_2(G) = G^{[2]}$\\
Clearly, $B$ corresponds to a quiver $Q$ on $3$ vertices, and $B^{[2]}=\mu_2(B).$ Notice that $k=2.$ Since the initial $C$-matrix $C$ and $G$-matrix $G$ are identity matrices, we have $C=G=I.$

Then $G^{[2]}=G(J_2 + [-\varepsilon_2(C)B]_+^{\bullet 2}).$
First, we check the sign. There exists $\varepsilon_2(C)= \pm 1$ such that $[-\varepsilon_2(C)C]_+^{\bullet 2}.$
In this case, $\varepsilon_2(C) = 1$ since this satisfies (\ref{Csigncoherence}) i.e. 
\begin{equation*}
[-C]_+^{\bullet 2}
= \left[\begin{smallmatrix}
-1 & 0 & 0 \\ 0 & -1 & 0 \\ 0 & 0 & -1
\end{smallmatrix}\right]_+^{\bullet 2}
= \left[\begin{smallmatrix}
0 & 0 & 0 \\ 0 & 0 & 0 \\ 0 & 0 & 0
\end{smallmatrix}\right]^{\bullet 2}
= \mathbf{0}.
\end{equation*}
Now we have
\begin{align*}
G^{[2]}&=G(J_2 + [-\varepsilon_2(C)B]_+^{\bullet 2})\\
&=G(J_2 + [-B]_+^{\bullet 2})\\
&=I\left(\left[\begin{smallmatrix}
1 & 0 & 0 \\ 0 & -1 & 0 \\ 0 & 0 & 1
\end{smallmatrix}\right] + [-B]_+^{\bullet 2}\right).\\
\end{align*}
Here $[-B]_+^{\bullet 2}$ is obtained by the following procedure:\\
$$[-B]_+^{\bullet 2} = \left[\begin{smallmatrix}
0 & 3 & -3 \\ -3 & 0 & 3 \\ 3 & -3 & 0
\end{smallmatrix}\right]_+^{\bullet 2}
= \left[\begin{smallmatrix}
0 & 3 & 0 \\ 0 & 0 & 3 \\ 3 & 0 & 0
\end{smallmatrix}\right]^{\bullet 2}
= \left[\begin{smallmatrix}
0 & 3 & 0 \\ 0 & 0 & 0 \\ 0 & 0 & 0
\end{smallmatrix}\right].$$
So we have
\begin{align*}
G^{[2]}&=I\left(\left[\begin{smallmatrix}
1 & 0 & 0 \\ 0 & -1 & 0 \\ 0 & 0 & 1
\end{smallmatrix}\right] + \left[\begin{smallmatrix}
0 & 3 & 0 \\ 0 & 0 & 0 \\ 0 & 0 & 0
\end{smallmatrix}\right]\right)
=\left[\begin{smallmatrix}
1 & 3 & 0 \\ 0 & -1 & 0 \\ 0 & 0 & 1
\end{smallmatrix}\right].
\end{align*}

\noindent\textbf{ii)} $\mu_3(G^{[2]}) = G^{[2,3]}$\\
Suppose we know that $B^{[2]}=\left[
\begin{smallmatrix}
0 & 3 & -6 \\ -3 & 0 & 3 \\ 6 & -3 & 0
\end{smallmatrix}\right]$, and now $B^{[2,3]}=\mu_3(B^{[2]}).$ Notice that $k=3.$ Also, we have $C^{[2]}=\left[\begin{smallmatrix}
1 & 0 & 0 \\ 3 & -1 & 0 \\ 0 & 0 & 1
\end{smallmatrix}\right],$ and $G^{[2]}=\left[\begin{smallmatrix}
1 & 3 & 0 \\ 0 & -1 & 0 \\ 0 & 0 & 1
\end{smallmatrix}\right]$ from the previous case.

Then $G^{[2,3]}=G^{[2]}(J_3 + [-\varepsilon_3(C^{[2]})B^{[2]}]_+^{\bullet 3}).$
First, we check the sign. There exists $\varepsilon_3(C^{[2]})= \pm 1$ such that $[-\varepsilon_3(C^{[2]})C^{[2]}]_+^{\bullet 3}.$
In this case, $\varepsilon_3(C^{[2]}) = 1$ since this satisfies (\ref{Csigncoherence}) i.e. 
\begin{equation*}
[-C^{[2]}]_+^{\bullet 3}
= \left[\begin{smallmatrix}
-1 & -3 & 0 \\ 0 & 1 & 0 \\ 0 & 0 & -1
\end{smallmatrix}\right]_+^{\bullet 3}
= \left[\begin{smallmatrix}
0 & 0 & 0 \\ 0 & 1 & 0 \\ 0 & 0 & 0
\end{smallmatrix}\right]^{\bullet 3}
= \mathbf{0}.
\end{equation*}
Now we have
\begin{align*}
G^{[2,3]}&=G^{[2]}(J_3 + [-\varepsilon_3(C^{[2]})B^{[2]}]_+^{\bullet 3})=G^{[2]}(J_3 + [-B^{[2]}]_+^{\bullet 3})= \left[\begin{smallmatrix}
1 & 3 & 0 \\ 0 & -1 & 0 \\ 0 & 0 & 1
\end{smallmatrix}\right]
\left(\left[\begin{smallmatrix}
1 & 0 & 0 \\ 0 & 1 & 0 \\ 0 & 0 & -1
\end{smallmatrix}\right] + [-B^{[2]}]_+^{\bullet 3}\right).
\end{align*}
Here $[-B^{[2]}]_+^{\bullet 3}$ is obtained by the following procedure:
$$[-B^{[2]}]_+^{\bullet 3} = \left[\begin{smallmatrix}
0 & -3 & 6 \\ 3 & 0 & -3 \\ -6 & 3 & 0
\end{smallmatrix}\right]_+^{\bullet 3}
= \left[\begin{smallmatrix}
0 & 0 & 6 \\ 3 & 0 & 0 \\ 0 & 3 & 0
\end{smallmatrix}\right]^{\bullet 3}
= \left[\begin{smallmatrix}
0 & 0 & 6 \\ 0 & 0 & 0 \\ 0 & 0 & 0
\end{smallmatrix}\right].$$
So we have
\begin{align*}
G^{[2,3]}&=\left[\begin{smallmatrix}
1 & 3 & 0 \\ 0 & -1 & 0 \\ 0 & 0 & 1
\end{smallmatrix}\right]
\left(\left[\begin{smallmatrix}
1 & 0 & 0 \\ 0 & 1 & 0 \\ 0 & 0 & -1
\end{smallmatrix}\right] + \left[\begin{smallmatrix}
0 & 0 & 6 \\ 0 & 0 & 0 \\ 0 & 0 & 0
\end{smallmatrix}\right]\right)
=\left[\begin{smallmatrix}
1 & 3 & 6 \\ 0 & -1 & 0 \\ 0 & 0 & -1
\end{smallmatrix}\right].\\
\end{align*}

\noindent\textbf{iii)} $\mu_1(G^{[2,3]}) = G^{[2,3,1]}$\\
Suppose we know that $B^{[2,3]}=\left[
\begin{smallmatrix}
0 & -15 & 6 \\ 15 & 0 & -3 \\ -6 & 3 & 0
\end{smallmatrix}\right]$, and now $B^{[2,3,1]}=\mu_1(B^{[2,3]}).$ Notice that $k=3.$ Also, we have $C^{[2,3]}=\left[\begin{smallmatrix}
1 & 0 & 0 \\ 3 & -1 & 0 \\ 6 & 0 & -1
\end{smallmatrix}\right],$ and $G^{[2,3]}=\left[\begin{smallmatrix}
1 & 3 & 6 \\ 0 & -1 & 0 \\ 0 & 0 & -1
\end{smallmatrix}\right]$ from the previous case.

Then, $G^{[2,3,1]}=G^{[2,3]}(J_1 + [-\varepsilon_1(C^{[2,3]})B^{[2,3]}]_+^{\bullet 1}).$
First, check the sign. There exists $\varepsilon_1(C^{[2,3]})= \pm 1$ such that $[-\varepsilon_1(C^{[2,3]})C^{[2,3]}]_+^{\bullet 1}.$
In this case, $\varepsilon_1(C^{[2,3]}) = 1$ since this satisfies (\ref{Csigncoherence}) i.e. 

\begin{equation*}
[-C^{[2,3]}]_+^{\bullet 1}
= \left[\begin{smallmatrix}
-1 & 0 & 0 \\ -3 & 1 & 0 \\ -6 & 0 & 1
\end{smallmatrix}\right]_+^{\bullet 1}
= \left[\begin{smallmatrix}
0 & 0 & 0 \\ 0 & 1 & 0 \\ 0 & 0 & 1
\end{smallmatrix}\right]^{\bullet 1}
= \mathbf{0}.
\end{equation*}

Now, we have
\begin{align*}
G^{[2,3,1]}&=G^{[2,3]}(J_1 + [-\varepsilon_1(C^{[2,3]})B^{[2,3]}]_+^{\bullet 1})\\
&=G^{[2,3]}(J_1 + [-B^{[2,3]}]_+^{\bullet 1})\\
&= \left[\begin{smallmatrix}
1 & 3 & 6 \\ 0 & -1 & 0 \\ 0 & 0 & -1
\end{smallmatrix}\right]
\left(\left[\begin{smallmatrix}
-1 & 0 & 0 \\ 0 & 1 & 0 \\ 0 & 0 & 1
\end{smallmatrix}\right] + [-B^{[2,3]}]_+^{\bullet 1}\right).\\
\end{align*}
Here $[-B^{[2,3]}]_+^{\bullet 1}$ is obtained by the following procedure:\\
$$[-B^{[2,3]}]_+^{\bullet 1} = \left[\begin{smallmatrix}
0 & 15 & -6 \\ -15 & 0 & 3 \\ 6 & -3 & 0
\end{smallmatrix}\right]_+^{\bullet 1}
= \left[\begin{smallmatrix}
0 & 15 & 0 \\ 0 & 0 & 3 \\ 6 & 0 & 0
\end{smallmatrix}\right]^{\bullet 1}
= \left[\begin{smallmatrix}
0 & 0 & 0 \\ 0 & 0 & 0 \\ 6 & 0 & 0
\end{smallmatrix}\right].$$

So we have
\begin{align*}
G^{[2,3,1]}&=\left[\begin{smallmatrix}
1 & 3 & 6 \\ 0 & -1 & 0 \\ 0 & 0 & -1
\end{smallmatrix}\right]
\left(\left[\begin{smallmatrix}
-1 & 0 & 0 \\ 0 & 1 & 0 \\ 0 & 0 & 1
\end{smallmatrix}\right] + \left[\begin{smallmatrix}
0 & 0 & 0 \\ 0 & 0 & 0 \\ 6 & 0 & 0
\end{smallmatrix}\right]\right)
=\left[\begin{smallmatrix}
35 & 3 & 6 \\ 0 & -1 & 0 \\ -6 & 0 & -1
\end{smallmatrix}\right].\\
\end{align*}

\section{Rank 3 mutation-cyclic quivers}

\begin{Def} (Mutation-cyclic)
We say that a quiver $Q$ is \emph{mutation-cyclic} if for any sequence of mutations $\bw$ applied to $Q$, the resulting quiver $Q^{\bw}$ is cyclic. That is, no matter how many times or in which order mutations are performed on $Q,$ the quiver retains the property of being cyclic. Formally, let $\mu_k$ denote the mutation at vertex $k.$ Then $Q$ is mutation cyclic if for any finite sequence of mutations $\mu_{i_1}, \mu_{i_2}, \ldots, \mu_{i_t}$ applied to $Q,$ the resulting quiver $\mu_{i_t} \cdots \mu_{i_2} \mu_{i_1}(Q)$ is cyclic.
\end{Def}


In order to state our theorem, we fix $Q$ as a rank 3 mutation-cyclic quiver. 

\begin{Def}\label{Amatdefn} ($A$-matrix)\\
Let $B$ be a $3\times 3$ skew-symmetric matrix. 
\emph{The pseudo-Cartan companion}, or $A$-matrix, of $B$ is the $3\times 3$ symmetric matrix $A=[a_{ij}]$ such that $a_{ii} = 2$ for all $i\in\{1,2,3\}$, and $a_{ij} = |b_{ij}|$ for all $i\neq j.$
\end{Def}

\begin{remark}
In \cite{seven2015cluster}, it is shown that for skew-symmetric cluster algebras, the c-vectors associated with an acyclic seed define a quasi-Cartan companion\footnote{The key difference between the pseudo-Cartan companion and the quasi-Cartan companion lies in the definition of the off-diagonal entries. In the quasi-Cartan companion, the off-diagonal entries $a_{ij}$ satisfy the condition $|a_{ij}| = |b_{ij}|$ for all $i \neq j$ while in the pseudo-Cartan companion, they are given by $a_{ij} = |b_{ij}|$.}. Here, we needed to define the pseudo-Cartan companion, denoted by $A$, as a fundamental tool to prove Theorem \ref{thm-rank3}.
\end{remark}

\begin{Def} ($A$-mutation)\\ \label{Adef} 
Let $B$ be a $3\times 3$ exchange matrix corresponding to rank 3 mutation cyclic quiver. 
For each mutation sequence $\bw=[i_1, i_2, \dots , i_t]$, we inductively define the $A$-matrix of $B$ as follows.

Assume $A$ is a pseudo-Cartan companion of the initial $B.$ Let $\bw=[i_1, i_2, \dots , i_t]$ be a mutation sequence. Suppose that $A^\bw$ is defined. We seek to define $A^{\bw[k]}$ such that $\mu_k(A^{\bw}) = A^{\bw[k]}.$

For each $l\in\{1,2,3\}\setminus\{k\}$, we define a $3\times 3$ matrix $D^{(k,l)}=\left[d_{rs}^{(k,l)}\right]$ as follows: 

\begin{equation*} 
d_{rs}^{(k,l)} = 
\begin{cases} a_{kl}^{\bw} & \text{if  $r=k, s=l$}, 
\\ -1 & \text{if $r=s=l$}, 
\\ 1 & \text{if $r=s\neq l$},
\\ 0 & \text{else}.
\end{cases}
\end{equation*}

These matrices \( D^{(k,l)} \) are used in the mutation formula to update \( A^{\bw} \) to \( A^{\bw[k]} \). Specifically, after mutating at \( k \), we have:

\[ A^{\bw[k]} = (D^{(k,l)})^T A^{\bw} D^{(k,l)}. \]

\end{Def}

We can now state our main results. 
When a rank 3 mutation cyclic quiver is given by $Q,$ let a reduced mutation sequence be $\bw=[i_1, i_2, \dots , i_t]$ with $i_r \in \{1,2,3\}.$ 
For simplicity, let $\mathbf{i}=i_1,$ the first mutation in a mutation sequence. Then after mutating at vertex $\mathbf{i}$, we will have a corresponding $B^{[\mathbf{i}]}=[b_{ij}^{[\mathbf{i}]}].$ 
\begin{theorem} \label{thm-rank3}
We have $A^\bw = (G^{\bw[m]})^T A^{[\mathbf{i}]} (G^{\bw[m]})$ for any reduced mutation sequence $\bw$ with length of at least 1, starting with $\mathbf{i}.$
\end{theorem}

\begin{corollary} \label{cor-rank3}
     Define a quadratic form 
     \[ q(x_1,x_2,x_3) =\sum_{i=1}^3 x_i^2 + \sum_{1\leq i<j\leq 3} b_{ij}^{[\mathbf{i}]} x_i x_j,\]
     where $b_{ij}^{[\mathbf{i}]} > 0$ is the number of arrow between $i$ and $j$ in $\mu_\mathbf{i}(Q).$
 Then for any reduced mutation sequence $\bw = [\mathbf{i}, i_2, \dots , i_t]$ and for any $i \in \{1,2,3\}$, the  column vector $\mathbf{g}_i^{\bw}$ in the corresponding $G$-matrix  is a solution to the following equation :
\begin{equation}
q(\mathbf{g}_{i}^\bw) = q\begin{pmatrix}
\begin{bmatrix}x_1 \\ x_2 \\ x_3 \end{bmatrix}\end{pmatrix} \\
= x_1^2+x_2^2+x_3^2 + |b^{[\mathbf{i}]}_{12}|x_1x_2 + |b^{[\mathbf{i}]}_{23}|x_2x_3 + |b^{[\mathbf{i}]}_{31}|x_3x_1 = 1.
\end{equation}
\end{corollary}

\subsection{Examples of Theorem \ref{thm-rank3} and Corollary \ref{cor-rank3}}

\begin{example}
Here is a specific example of the identity $A^\bw = (G^{\bw[m]})^T A^{[\mathbf{i}]} (G^{\bw[m]}).$ 

Start with a case when $\bw = [2,3],$ and $m=1.$ Then $\mathbf{i} = 2$ in this case. We are going to show the identity $A^{[2,3]} = (G^{[2,3,1]})^T A^{[2]} (G^{[2,3,1]})$ as an example. 

If we take \ref{rank3initialquiver} as an our initial quiver, $\mu_2(Q) =  Q^{[2]} $ can be depicted as
\begin{equation}
\begin{tikzpicture}[->,>=stealth',shorten >=1pt,auto,node distance=1.5cm,thick,main node/.style={circle,draw}]
  \node[main node] (2) {2};
  \node[main node] (1) [below left of=2] {1};
  \node[main node] (3) [below right of=2] {3};

  \path
    (2) edge [left] node[auto, midway, above] {3} (1)  
    (3) edge [left] node[auto, midway, above] {3} (2)
    (1) edge [right] node[auto, midway, above] {3} (3);
\end{tikzpicture}
\hspace{1cm}
\xrightarrow{\text{$\mu_2$}}
\hspace{1cm}
\begin{tikzpicture}[->,>=stealth',shorten >=1pt,auto,node distance=1.5cm,thick,main node/.style={circle,draw}]
  \node[main node] (2) {2};
  \node[main node] (1) [below left of=2] {1};
  \node[main node] (3) [below right of=2] {3};

  \path
    (1) edge [left] node[auto, midway, above] {3} (2)  
    (2) edge [left] node[auto, midway, above] {3} (3)
    (3) edge [right] node[auto, midway, above] {6} (1);
\end{tikzpicture}
\end{equation}

Observe that $\mu_2(Q)$ has 3 arrows from vertex 1 to vertex 2, 3 arrows from vertex 2 to vertex 3, 6 arrows from vertex 3 to vertex 1, i.e., $b_{12}^{[2]} = 3, b_{23}^{[2]} = 3, b_{31}^{[2]} = 6.$ So, by the Definition \ref{Amatdefn}, $A^{[2]} = \left[\begin{smallmatrix}
2 & 3 & 6 \\3 & 2 & 3 \\ 6 & 3 & 2
\end{smallmatrix}\right].$

From subsection \ref{example2.1}, we have the $G$-matrix for sequence $\bw = [2,3,1],$ that is $G^{[2,3,1]}=\left[\begin{smallmatrix}
35 & 3 & 6 \\ 0 & -1 & 0 \\ -6 & 0 & -1
\end{smallmatrix}\right].$
Thus, we can easily get that $(G^{[2,3,1]})^T A^{[2]} (G^{[2,3,1]})= \left[\begin{smallmatrix}
35 & 3 & 6 \\ 0 & -1 & 0 \\ -6 & 0 & -1
\end{smallmatrix}\right]^T
\left[\begin{smallmatrix}
2 & 3 & 6 \\3 & 2 & 3 \\ 6 & 3 & 2
\end{smallmatrix}\right]
\left[\begin{smallmatrix}
35 & 3 & 6 \\ 0 & -1 & 0 \\ -6 & 0 & -1
\end{smallmatrix}\right]
= \left[\begin{smallmatrix}
2 & 15 & 6 \\ 15 & 2 & 3 \\ 6 & 3 & 2
\end{smallmatrix}\right].$\\

On the other hand, $A^{[2,3]}$ can be obtained by using either $D^{(3,1)} = \left[\begin{smallmatrix}
-1 & 0 & 0 \\ 0 & 1 & 0 \\ 6 & 0 & 1
\end{smallmatrix}\right]$ 
or
$D^{(3,2)} = \left[\begin{smallmatrix}
1 & 0 & 0 \\ 0 & -1 & 0 \\ 0 & 3 & 1
\end{smallmatrix}\right] : $

\begin{equation*}
A^{[2,3]} = (D^{(3,1)})^T A^{[2]} D^{(3,1)} = \left[\begin{smallmatrix}
-1 & 0 & 0 \\ 0 & 1 & 0 \\ 6 & 0 & 1
\end{smallmatrix}\right]^T
\left[\begin{smallmatrix}
2 & 3 & 6 \\3 & 2 & 3 \\ 6 & 3 & 2
\end{smallmatrix}\right]
\left[\begin{smallmatrix}
-1 & 0 & 0 \\ 0 & 1 & 0 \\ 6 & 0 & 1
\end{smallmatrix}\right] =
\left[\begin{smallmatrix}
2 & 15 & 6 \\ 15 & 2 & 3 \\ 6 & 3 & 2
\end{smallmatrix}\right]
\end{equation*}

\begin{equation*}
A^{[2,3]} = (D^{(3,2)})^T A^{[2]} D^{(3,2)} = \left[\begin{smallmatrix}
1 & 0 & 0 \\ 0 & -1 & 0 \\ 0 & 3 & 1
\end{smallmatrix}\right]^T
\left[\begin{smallmatrix}
2 & 3 & 6 \\3 & 2 & 3 \\ 6 & 3 & 2
\end{smallmatrix}\right]
\left[\begin{smallmatrix}
1 & 0 & 0 \\ 0 & -1 & 0 \\ 0 & 3 & 1
\end{smallmatrix}\right] = \left[\begin{smallmatrix}
2 & 15 & 6 \\ 15 & 2 & 3 \\ 6 & 3 & 2
\end{smallmatrix}\right] 
\end{equation*}

Hence,
\begin{equation*}
A^{[2,3]} = (G^{[2,3,1]})^T A^{[2]} (G^{[2,3,1]}).     
\end{equation*}
Furthermore,
\begin{align*}
A^{[2,3,1]}
&= (G^{[2,3,1,2]})^T A^{[2]} (G^{[2,3,1,2]}) = \left[\begin{smallmatrix}
35 & 522 & 6 \\ 0 & 1 & 0 \\ -6 & -90 & -1
\end{smallmatrix}\right]^T
\left[\begin{smallmatrix}
2 & 3 & 6 \\3 & 2 & 3 \\ 6 & 3 & 2
\end{smallmatrix}\right]
\left[\begin{smallmatrix}
35 & 522 & 6 \\ 0 & 1 & 0 \\ -6 & -90 & -1
\end{smallmatrix}\right]
= \left[\begin{smallmatrix}
2 & 15 & 6 \\ 15 & 2 & 87 \\ 6 & 87 & 2
\end{smallmatrix}\right] \\
A^{[2,3,1,2]} &= (G^{[2,3,1,2,1]})^T A^{[2]} (G^{[2,3,1,2,1]}) = \left[\begin{smallmatrix}
7795 & 522 & 6 \\ 15 & 1 & 0 \\ -1344 & -90 & -1
\end{smallmatrix}\right]^T
\left[\begin{smallmatrix}
2 & 3 & 6 \\3 & 2 & 3 \\ 6 & 3 & 2
\end{smallmatrix}\right]
\left[\begin{smallmatrix}
7795 & 522 & 6 \\ 15 & 1 & 0 \\ -1344 & -90 & -1
\end{smallmatrix}\right]
=\left[\begin{smallmatrix}
2 & 15 & 1299 \\ 15 & 2 & 87 \\ 1299 & 87 & 2
\end{smallmatrix}\right]\\
\end{align*}
\end{example}

\begin{example}
Again, in the context of subsection \ref{example2.1}, we have the following $G$-matrices for each sequence respectively : 
$G=\left[\begin{smallmatrix}
1 & 0 & 0 \\ 0 & 1 & 0 \\ 0 & 0 & 1
\end{smallmatrix}\right],
G^{[2]}=\left[\begin{smallmatrix}
1 & 3 & 0 \\ 0 & -1 & 0 \\ 0 & 0 & 1
\end{smallmatrix}\right], G^{[2,3]}=\left[\begin{smallmatrix}
1 & 3 & 6 \\ 0 & -1 & 0 \\ 0 & 0 & -1
\end{smallmatrix}\right],$ and $G^{[2,3,1]}=\left[\begin{smallmatrix}
35 & 3 & 6 \\ 0 & -1 & 0 \\ -6 & 0 & -1
\end{smallmatrix}\right].$
See that each of those mutation sequences has $\mathbf{i}=2$ as the first mutation in the mutation sequence. 

As we have $b_{12}^{[2]} = 3, b_{23}^{[2]} = 3, b_{31}^{[2]} = 6,$ each $\mathbf{g}_i^{\bw} = \begin{bmatrix}x_1 \\ x_2 \\ x_3 \end{bmatrix}$ from $G^{\bw}$-matrix is a solution to the equation $$q(\mathbf{g}_i^{\bw}) = q\begin{pmatrix}
\begin{bmatrix}x_1 \\ x_2 \\ x_3 \end{bmatrix}\end{pmatrix}
= x_1^2+x_2^2+x_3^2+3x_1x_2+3x_2x_3+6x_3 x_1 = 1,$$ i.e. 
\begin{equation*}
q(1,0,0)=q(0,1,0)=q(0,0,1)=q(3,-1,0)=q(6,0,-1)=q(35,0,-6) = 1. 
\end{equation*}
\end{example}

\section{Proof of Theorem \ref{thm-rank3}}

We will prove that $A^\bw = (G^{\bw[m]})^T A^{[\mathbf{i}]} (G^{\bw[m]})$ for any reduced mutation sequence $\bw$ with a length of at least 1, starting with $\mathbf{i},$ by induction.
Let us note that for $\mu_m(G^\bw, A^\bw) = (G^{\bw[m]}, A^{\bw[m]}).$

Our proof consists of two parts, which together imply Theorem \ref{thm-rank3}.
\begin{itemize}
    \item In subsection \ref{sub4.1}, we will show the base step of our induction. Specifically, we will porve the identity $A^\bw = (G^{\bw[m]})^T A^{[\mathbf{i}]} (G^{\bw[m]})$ for $\bw = [\mathbf{i}]$ with $\mathbf{i} \in \{1,2,3\},$ $m \in \{1,2,3\}\setminus\{\mathbf{i}\}$

    \item In subsection \ref{sub4.2}, we will show that for $|\bw| \geq 1,$ starting with $\mathbf{i},$ if the identity $A^\bw = (G^{\bw[m]})^T A^{[\mathbf{i}]} (G^{\bw[m]})$ holds, then the identity $A^{\bw[m]} = (G^{\bw[m,n]})^T A^{[\mathbf{i}]} (G^{\bw[m,n]})$ also holds for $m \neq n.$
\end{itemize}

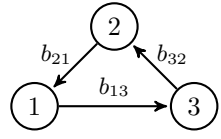
\begin{equation}\label{thmrank3initialquiver}
\begin{tikzpicture}[->,>=stealth',shorten >=1pt,auto,node distance=1.5cm,thick,main node/.style={circle,draw}]
  \node[main node] (2) {2};
  \node[main node] (1) [below left of=2] {1};
  \node[main node] (3) [below right of=2] {3};

  \path
    (2) edge [left] node[auto, pos=0.9, above, yshift=4pt] {\small $b_{21}$} (1)  
    (3) edge [left] node[auto, pos=0.1, above, yshift=4pt] {\small $b_{32}$} (2)
    (1) edge [right] node[auto, midway, above] {\small $b_{13}$} (3);
\end{tikzpicture}
\end{equation}

Let $B =\left[\begin{smallmatrix}
0 & -b_{12} & b_{13} \\ b_{12} & 0 & -b_{23} \\ -b_{13} & b_{23} & 0    
\end{smallmatrix}\right]$ be an initial exchange matrix corresponding to rank 3 mutation cycllic quiver (\ref{thmrank3initialquiver}), where $b_{12},b_{13},b_{23} > 0.$

\subsection{Base Step}\label{sub4.1}

Let $\bw = [\mathbf{i}]$ be a mutation sequence and let us note that for 
$\mu_n(G^{\bw[m]}, A^{\bw[m]}) = (G^{\bw[m,n]},A^{\bw[m,n]}).$
We will show the identity $A^{[\mathbf{i}]} = (G^{[\mathbf{i},m]})^T A^{[\mathbf{i}]} (G^{[\mathbf{i},m]})$ for the cases where ($\mathbf{i},m) = (1,2)$ and $(\mathbf{i},m) = (1,3)$. The remaining cases, $(\mathbf{i},m) = (2,1),(2,3),(3,1),$ and $(3,2)$, can be verified through analogous steps.\\

Recall that $G=I,$ $G^{[\mathbf{i}]} = G (J_{\mathbf{i}} + [-\varepsilon_{\mathbf{i}}(C)B]_+^{\bullet \mathbf{i}}),$ and 
$G^{[\mathbf{i},m]} = G^{[\mathbf{i}]} (J_{m} + [-\varepsilon_{m}(C^{[\mathbf{i}]})B^{[\mathbf{i}]}]_+^{\bullet m})$
by (\ref{CG-mutation}). In what follows, LHS (resp. RHS) stands for the left hand side (resp. right hand side).
We consider the case of $\mathbf{i}=1, m=2$, as the other cases are similar. 

\noindent\fbox{LHS} $A^{[1]}$\\
As $A = \left[\begin{smallmatrix}
2 & a_{12} & a_{13} \\ a_{12} & 2 & a_{23} \\ a_{13} & a_{23} & 2
\end{smallmatrix}\right],$ by the mutation of $A,$ we can easily get
$$A^{[1]} = \left[\begin{smallmatrix}
2 & a_{12} & a_{13} \\ a_{12} & 2 &  a_{12} a_{13}-a_{23} \\ a_{13}& a_{12}a_{13}-a_{23}& 2
\end{smallmatrix}\right].$$\\
\fbox{RHS} $(G^{[1,2]})^T A^{[1]} (G^{[1,2]})$\\
As $G = I$ and $\varepsilon_{1}(C) = 1,$ 
\begin{align*}
G^{[1]} & = G(J_{1} + [-\varepsilon_{1}(C)B]_+^{\bullet 1})\\
& = G(J_1 + [-B]_+^{\bullet 1})\\
& = \left[\begin{smallmatrix}
1 & 0 & 0 \\ 0 & 1 & 0 \\ 0 & 0 & 1
\end{smallmatrix}\right]
\left(\left[\begin{smallmatrix}
-1 & 0 & 0 \\ 0 & 1 & 0 \\ 0 & 0 & 1
\end{smallmatrix}\right] + [-B]_+^{\bullet 1}\right)\\
\end{align*}
Here, $[-B]_+^{\bullet 1}$ is obtained by the following procedure:\\
$$[-B]_+^{\bullet 1} = \left[\begin{smallmatrix}
0 & b_{12} & -b_{13} \\ -b_{12} & 0 & b_{23} \\ b_{13} & -b_{23} & 0    
\end{smallmatrix}\right]_+^{\bullet 1}
= \left[\begin{smallmatrix}
0 & b_{12} & 0 \\ 0 & 0 & b_{23} \\ b_{13} & 0 & 0    
\end{smallmatrix}\right]^{\bullet 1}
= \left[\begin{smallmatrix}
0 & 0 & 0 \\ 0 & 0 & 0 \\ b_{13} & 0 & 0
\end{smallmatrix}\right]$$
So,
\begin{equation*}
G^{[1]} = \left[\begin{smallmatrix}
1 & 0 & 0 \\ 0 & 1 & 0 \\ 0 & 0 & 1
\end{smallmatrix}\right]
\left(\left[\begin{smallmatrix}
-1 & 0 & 0 \\ 0 & 1 & 0 \\ 0 & 0 & 1
\end{smallmatrix}\right] + \left[\begin{smallmatrix}
0 & 0 & 0 \\ 0 & 0 & 0 \\ b_{13} & 0 & 0
\end{smallmatrix}\right]\right)
=\left[\begin{smallmatrix}
1 & 0 & 0 \\ 0 & 1 & 0 \\ 0 & 0 & 1
\end{smallmatrix}\right]\left[\begin{smallmatrix}
-1 & 0 & 0 \\ 0 & 1 & 0 \\ b_{13} & 0 & 1
\end{smallmatrix}\right]
=\left[\begin{smallmatrix}
-1 & 0 & 0 \\ 0 & 1 & 0 \\ b_{13} & 0 & 1
\end{smallmatrix}\right]
\end{equation*}

Now, as we have $G^{[1]} = \left[\begin{smallmatrix}
-1 & 0 & 0 \\ 0 & 1 & 0 \\ b_{13} & 0 & 1
\end{smallmatrix}\right],$ $\varepsilon_{2}(C^{[1]}) = 1,$ and $B^{[1]} =\left[\begin{smallmatrix}
0 & b_{12} & -b_{13} \\ -b_{12} & 0 & b_{12}b_{13}-b_{23} \\ b_{13} & -(b_{12}b_{13}-b_{23}) & 0    
\end{smallmatrix}\right]$ by the definition of mutation of $B$, 
\begin{align*}
G^{[1,2]} & = G^{[1]} (J_{2} + [-\varepsilon_{2}(C^{[1]})B^{[1]}]_+^{\bullet 2})\\
& = G^{[1]}(J_2 + [-B^{[1]}]_+^{\bullet 2})\\
& = \left[\begin{smallmatrix}
-1 & 0 & 0 \\ 0 & 1 & 0 \\ b_{13} & 0 & 1
\end{smallmatrix}\right]
\left(\left[\begin{smallmatrix}
1 & 0 & 0 \\ 0 & -1 & 0 \\ 0 & 0 & 1
\end{smallmatrix}\right] + [-B^{[1]}]_+^{\bullet 2}\right)
\end{align*}
Here, $[-B^{[1]}]_+^{\bullet 2}$ is obtained by the following procedure:
$$[-B^{[\,]}]_+^{\bullet 2}
= \left[\begin{smallmatrix}
0 & -b_{12} & b_{13} \\ b_{12} & 0 & -(b_{12}b_{13}-b_{23}) \\ -b_{13} & b_{12}b_{13}-b_{23} & 0    
\end{smallmatrix}\right]_+^{\bullet 2}
= \left[\begin{smallmatrix}
0 & 0 & b_{13} \\ b_{1,2} & 0 & 0 \\ 0 & b_{12}b_{13}-b_{23} & 0      
\end{smallmatrix}\right]^{\bullet 2}
= \left[\begin{smallmatrix}
0 & 0 & 0 \\ 0 & 0 & 0 \\ 0 & b_{12}b_{13}-b_{23} & 0
\end{smallmatrix}\right]$$

So,
\begin{equation*}
G^{[1,2]} = \left[\begin{smallmatrix}
-1 & 0 & 0 \\ 0 & 1 & 0 \\ b_{13} & 0 & 1
\end{smallmatrix}\right]
\left(\left[\begin{smallmatrix}
1 & 0 & 0 \\ 0 & -1 & 0 \\ 0 & 0 & 1
\end{smallmatrix}\right] + \left[\begin{smallmatrix}
0 & 0 & 0 \\ 0 & 0 & 0 \\ 0 & b_{12}b_{13}-b_{23} & 0
\end{smallmatrix}\right]\right)
= \left[\begin{smallmatrix}
-1 & 0 & 0 \\ 0 & 1 & 0 \\ b_{13} & 0 & 1
\end{smallmatrix}\right] \left[\begin{smallmatrix}
1 & 0 & 0 \\ 0 & -1 & 0 \\ 0 & b_{12}b_{13}-b_{23} & 1
\end{smallmatrix}\right]
= \left[\begin{smallmatrix}
-1 & 0 & 0 \\ 0 & -1 & 0 \\ b_{13} & b_{12}b_{13}-b_{23} & 1
\end{smallmatrix}\right]
\end{equation*}
Then, since $b_{12}, b_{13}, b_{23} > 0$, we have $b_{12} = a_{12}$, $b_{13} = a_{13}$, and $b_{23} = a_{23}$. And thus, \begin{align*}
&(G^{[1,2]})^T A^{[1]}\\
&= \left[\begin{matrix}
-1 & 0 & b_{13} \\ 0 & -1 & b_{12}b_{13}-b_{23} \\ 0 & 0 & 1\end{matrix}\right]
\left[\begin{matrix}
2 & a_{12} & a_{13} \\ a_{12} & 2 & a_{12}a_{13}-a_{23} \\ a_{13} & a_{12}a_{13}-a_{23} & 2
\end{matrix}\right]\\
&= \left[\begin{matrix}
-2+b_{13}a_{13} & -a_{12}+b_{13}(a_{12}a_{13}-a_{23}) & -a_{13}+2b_{13} \\ -a_{12}+a_{13}(b_{12}b_{13}-b_{23}) & -2+(b_{12}b_{13}-b_{23})(a_{12}a_{13}-a_{23}) & -(a_{12}a_{13}-a_{23})+2(b_{12}b_{13}-b_{23}) \\ a_{13} & a_{12}a_{13}-a_{23} & 2
\end{matrix}\right]\\
&=\left[\begin{matrix}
-2+(a_{13})^2 & -a_{12}+a_{12}(a_{13})^2-a_{13}a_{23} & a_{13} \\ -a_{12}+a_{12}(a_{13})^2-a_{13}a_{23} & -2+(a_{12}a_{13}-a_{23})^2 & a_{12}a_{13}-a_{23} \\ A_{13} & a_{12}a_{13}-a_{23} & 2
\end{matrix}\right]\\
\end{align*}  
So
\begin{align*}
&(G^{[1,2]})^T A^{[1]} G^{[1,2]}\\
&=\left[\begin{matrix}
-2+(a_{13})^2 & -a_{12}+a_{12}(a_{13})^2-a_{13}a_{23} & a_{13} \\ -a_{12}+a_{12}(a_{13})^2-a_{13}a_{23} & -2+(a_{12}a_{13}-a_{23})^2 & a_{12}a_{13}-a_{23} \\ a_{13} & a_{12}a_{13}-a_{23} & 2
\end{matrix}\right]
\left[\begin{matrix}
-1 & 0 & 0 \\ 0 & -1 & 0 \\ b_{13} & b_{12}b_{13}-b_{23} & 1
\end{matrix}\right]\\
&=\left[\begin{matrix}
2 & a_{12} & a_{13} \\ a_{12} & 2 & a_{12}a_{13}-a_{23} \\ a_{13} & a_{12}a_{13}-a_{23} & 2
\end{matrix}\right]\\
&=A^{[1]}\\
\end{align*}

\subsection{Inductive Step} \label{sub4.2}
Recall that
$$A^{\bw} \xrightarrow{\text{$\mu_{m}$}} A^{\bw[m]} = (D^{(m,l)})^T A^\bw D^{(m,l)} \quad\text{for}\quad l\in\{1,2,3\}\setminus\{m\}$$

Define $G^{\bw[m]}$ such that $\mu_m(G^\bw) = G^{\bw[m]}$ for a reduced mutation sequence $\bw=[i_1, i_2, \dots , i_t]$ with $i_r \in \{1,2,3\}.$ For simplicity, let $p=i_t,$ the last mutation in the sequence.

Rewrite $G$-matrix mutation notation (\ref{CG-mutation}) as
\begin{equation}\label{CG-mutation2}
G^{\bw[m]} = G^{\bw} (J_m + [-sgn(\mathbf{c}_m^{\bw}) b_{pm}^{\bw}]_+ E^{pm}),
\end{equation}
where $J_l$ is the diagonal matrix obtained from the identity matrix by replacing the $(l,l)-$entry with $-1$, and
where $E^{ij}$ is the matrix whose only non-zero entry is a 1 in the $i^{th}$ row and $j^{th}$ column.

If you mutate at $u \neq m$ consecutively, then we have
$$G^{\bw[m]} \xrightarrow{\text{$\mu_{u}$}} G^{\bw[m,u]} = G^{\bw[m]}(J_u + [-sgn(\mathbf{c}_u^{\bw[m]}) b_{mu}^{\bw[m]}]_+ E^{mu}).$$

Suppose that $A^\bw = (G^{\bw[m]})^T A^{[\mathbf{i}]} (G^{\bw[m]})$ holds for $|\bw| \geq 1.$ 
Keep in mind that
$G^{\bw[m,u]} = G^{\bw[m]}(J_u + [-sgn(\mathbf{c}_u^{\bw[m]}) b_{mu}^{\bw[m]}]_+ E^{mu})$ and $A^{\bw[m]} = (D^{(m,l)})^T A^{\bw} D^{(m,l)}.$ \\

Let $N = J_u + [-sgn(\mathbf{c}_u^{\bw[m]}) b_{mu}^{\bw[m]}]_+ E^{mu}.$ Then, we will derive that $A^{\bw[m]} = (G^{\bw[m,u]})^T A^{\mathbf{i}} (G^{\bw[m,u]})$ by showing that $D^{(m,l)} = N.$

On the other hand, \footnote{Note that since $B$ is skew-symmetric, this results in a similar phenomenon as described in \cite[Theorem 1.3]{seven2015cluster}, even though the assumptions are different.} $[-sgn(\mathbf{c}_u^{\bw[m]}) b_{mu}^{\bw[m]}]_+ \neq 0$ if and only if $-sgn(\mathbf{c}_u^{\bw[m]}) b_{mu}^{\bw[m]} > 0$ if and only if $sgn(\mathbf{c}_u^{\bw[m]}) \neq sgn(b_{mu}^{\bw[m]}).$ Then we can say that $sgn(b_{mu}^{\bw[m]})=-sgn(\mathbf{c}_u^{\bw[m]}).$

Thus, \begin{align*}
-sgn(\mathbf{c}_u^{\bw[m]})b_{mu}^{\bw[m]}
&= sgn(b_{mu}^{\bw[m]}) b_{mu}^{\bw[m]} \\
&= |b_{mu}^{\bw[m]}|
\end{align*}


Mutating at vertex $m$ again (undoing the previous mutation) does not affect the number of arrows from vertex $m$ to vertex $u$, as these arrows are adjacent to vertex $m$. Only the direction of the arrows changes. Therefore,
\begin{equation*}
|b_{mu}^{\bw[m]}| = |-b_{mu}^{\bw}| = a_{mu}^{\bw} 
\end{equation*}

As $u$ depends on $m$ since $u$ cannot be the same as $m,$ for $\bw[m,u] = [i_1, i_2, \dots, i_t, m, u]$, we will check the case where $m = 1$ and $u \in \{2,3\}$. 

\noindent\fbox{LHS} $A^{\bw} \xrightarrow{\text{$\mu_{1}$}} A^{\bw[1]} = (D^{(1,l)})^T A^\bw D^{(1,l)}$
As $l$ is not equal to $m$ by the Definition \ref{Adef}, $l$ could be either 2 or 3 in this case. Thus, we will have two types of $D$ matrices, which are $D^{(1,2)},D^{(1,3)}$.

\underline{Type 1:}
\begin{equation*} 
d_{rs}^{(1,2)} = 
\begin{cases} a_{12}^{\bw} & \text{if  $r=1, s=2$}, 
\\ -1 & \text{if $r=s=2$}, 
\\ 1 & \text{if $r=s\neq 2$},
\\ 0 & \text{else}
\end{cases}
\quad \text{so,} \quad
D^{(1,2)}= \left[
\begin{matrix}
1 & a_{12}^{\bw} & 0 \\ 0 & -1 & 0 \\ 0 & 0 & 1
\end{matrix}
\right].
\end{equation*}

\underline{Type 2:}
\begin{equation*} 
d_{rs}^{(1,3)} = 
\begin{cases} a_{13}^{\bw} & \text{if  $r=1, s=3$}, 
\\ -1 & \text{if $r=s=3$}, 
\\ 1 & \text{if $r=s\neq 3$},
\\ 0 & \text{else}
\end{cases}
\quad \text{so,} \quad
D^{(1,3)}= \left[
\begin{matrix}
1 & 0 & a_{13}^{\bw} \\ 0 & 1 & 0 \\ 0 & 0 & -1
\end{matrix}
\right].
\end{equation*}
\\
\fbox{RHS} $G^{\bw[1]} \xrightarrow{\text{$\mu_{u}$}} G^{\bw[1,u]} = G^{\bw[1]}(J_u + [-sgn(\mathbf{c}_u^{\bw[1]}) b_{1u}^{\bw[1]}]_+ E^{1u}).$\\

\noindent Case $\underline{(m,u) = (1,2) :}$\\
Let $N_1$ be $J_2 + [-sgn(\mathbf{c}_2^{\bw[1]}) b_{12}^{\bw[1]}]_+ E^{12}.$ Then,
\begin{align*}
N_1 
&= J_2 + [-sgn(\mathbf{c}_2^{\bw[1]}) b_{12}^{\bw[1]}]_+ E^{12}\\
&= J_2 + |b_{12}^{\bw[1]}| E^{12}\\
&= J_2 + |-b_{12}^{\bw}| E^{12}\\
&= J_2 + a_{12}^{\bw} E^{12}\\
&= \left[\begin{smallmatrix}
1 & 0 & 0 \\ 0 & -1 & 0 \\ 0 & 0 & 1
\end{smallmatrix}\right]+
\left( a_{12}^{\bw}
\left[\begin{smallmatrix}
0 & 1 & 0 \\ 0 & 0 & 0 \\ 0 & 0 & 0
\end{smallmatrix}\right] 
\right)= \left[\begin{smallmatrix}
1 & a_{12}^{\bw} & 0 \\ 0 & -1 & 0 \\ 0 & 0 & 1
\end{smallmatrix}\right]\\
\end{align*}

\noindent Case $\underline{(m,u) = (1,3) :}$\\
Let $N_2$ be $J_3 + [-sgn(\mathbf{c}_3^{\bw[1]}) b_{13}^{\bw[1]}]_+ E^{13}.$ Then, 
\begin{align*}
N_2
&= J_3 + [-sgn(\mathbf{c}_3^{\bw[1]}) b_{13}^{\bw[1]}]_+ E^{13}\\
&= J_3 + |b_{13}^{\bw[1]}| E^{13} = J_3 + |-b_{13}^{\bw}| E^{13}\\
&= J_3 + a_{13}^{\bw} E^{13}\\
&= \left[\begin{smallmatrix}
1 & 0 & 0 \\ 0 & 1 & 0 \\ 0 & 0 & -1
\end{smallmatrix}\right]+
\left( a_{13}^{\bw}
\left[\begin{smallmatrix}
0 & 0 & 1 \\ 0 & 0 & 0 \\ 0 & 0 & 0
\end{smallmatrix}\right] 
\right)= \left[\begin{smallmatrix}
1 & 0 & a_{13}^{\bw} \\ 0 & 1 & 0 \\ 0 & 0 & -1
\end{smallmatrix}\right].
\end{align*}
Observe that
$$
N_1 = \left[\begin{smallmatrix}
1 & a_{12}^{\bw} & 0 \\ 0 & -1 & 0 \\ 0 & 0 & 1
\end{smallmatrix}\right]
= D^{(1,2)}\quad \text{ and }
\quad
N_2 = \left[
\begin{smallmatrix}
1 & 0 & a_{13}^{\bw} \\ 0 & 1 & 0 \\ 0 & 0 & -1
\end{smallmatrix}
\right] 
= D^{(1,3)}.$$
Eventually, we will get $A^{\bw[1]}$ from both cases above : 
\begin{align*}
(D^{(1,2)})^T A^\bw D^{(1,2)}
&= \left[\begin{smallmatrix}
1 & 0 & 0 \\ a_{12}^{\bw} & -1 & 0 \\ 0 & 0 & 1
\end{smallmatrix}\right]
\left[\begin{smallmatrix}
2 & a_{12}^{\bw} & a_{1,3}^{\bw} \\ a_{12}^{\bw} & 2 & a_{23}^{\bw} \\ a_{13}^{\bw} & a_{23}^{\bw} & 2
\end{smallmatrix}\right]
\left[\begin{smallmatrix}
1 & a_{12}^{\bw} & 0 \\ 0 & -1 & 0 \\ 0 & 0 & 1
\end{smallmatrix}\right]\\
&= \left[\begin{smallmatrix}
2 & a_{12}^{\bw} & a_{1,3}^{\bw} \\ a_{12}^{\bw} & (a_{12}^{\bw})^2 -2 & a_{12}^{\bw}a_{13}^{\bw}-a_{23}^{\bw} \\ a_{13}^{\bw} & a_{23}^{\bw} & 2
\end{smallmatrix}\right]
\left[\begin{smallmatrix}
1 & a_{12}^{\bw} & 0 \\ 0 & -1 & 0 \\ 0 & 0 & 1
\end{smallmatrix}\right]\\
&= \left[\begin{smallmatrix}
2 & a_{12}^{\bw} & a_{13}^{\bw} \\ a_{12}^{\bw} & 2 & a_{12}^{\bw}a_{13}^{\bw}-a_{23}^{\bw} \\ a_{13}^{\bw} & a_{12}^{\bw}a_{13}^{\bw}-a_{23}^{\bw} & 2
\end{smallmatrix}\right]=A^{\bw[1]}.\\
\end{align*}
\begin{align*}
(D^{(1,3)})^T A^\bw D^{(1,3)}
&= \left[\begin{smallmatrix}
1 & 0 & 0 \\ 0 & 1 & 0 \\ a_{13}^{\bw} & 0 & -1
\end{smallmatrix}\right]
\left[\begin{smallmatrix}
2 & a_{12}^{\bw} & a_{13}^{\bw} \\ a_{12}^{\bw} & 2 & a_{23}^{\bw} \\ a_{13}^{\bw} & a_{23}^{\bw} & 2
\end{smallmatrix}\right]
\left[\begin{smallmatrix}
1 & 0 & a_{13}^{\bw} \\ 0 & 1 & 0 \\ 0 & 0 & -1
\end{smallmatrix}\right]
\\
&= \left[\begin{smallmatrix}
2 & a_{12}^{\bw} & a_{13}^{\bw} \\ a_{12}^{\bw} & 2 & a_{23}^{\bw} \\ a_{13}^{\bw} & a_{12}^{\bw}a_{13}^{\bw}-a_{23}^{\bw} & (a_{13}^{\bw})^2-2
\end{smallmatrix}\right]
\left[\begin{smallmatrix}
1 & 0 & a_{13}^{\bw} \\ 0 & 1 & 0 \\ 0 & 0 & -1
\end{smallmatrix}\right]\\
&= \left[\begin{smallmatrix}
2 & a_{12}^{\bw} & a_{13}^{\bw} \\ a_{12}^{\bw} & 2 & a_{12}^{\bw}a_{13}^{\bw}-a_{23}^{\bw} \\ a_{13}^{\bw} & a_{12}^{\bw}a_{13}^{\bw}-a_{23}^{\bw} & 2
\end{smallmatrix}\right]=A^{\bw[1]}.
\end{align*}  

Similarly, The remaining cases, when $m = 2$ and $m = 3$ can be verified through analogous steps. Therefore, $D^{(m,l)} = N.$\\

By the induction, \ref{sub4.1} and \ref{sub4.2}, we can conclude that $A^{\bw[m]} = (G^{\bw[m,u]})^T A^{\mathbf{i}} (G^{\bw[m,u]})$for any reduced mutation sequence $\bw$ with a length of at least 1, starting with $\mathbf{i}.$  

\bibliographystyle{abbrv}
\bibliography{references.bib} 
\end{document}